\documentclass[a4,11pt]{amsart}
\usepackage{epsfig}
\usepackage{graphicx}
\usepackage{amssymb,amsthm, amsmath}
\def\N{\mathbb N}
\def\Z{\mathbb Z}
\def\Q{\mathbb Q}
\def\R{\mathbb R}

\def\A{\mathcal{A}}
\def\T{\mathcal T}
\def\L{\mathcal L}
\def\F{\mathcal F}
\def\S{\mathcal S}
\def\eps{\varepsilon}
\def\v{\mathbf{v}}
\newcommand{\supp}{\mbox{\rm supp}}

\newcommand{\LimL}{\underline{\mathrm{Lim}}}
\newcommand{\LimU}{\overline{\mathrm{Lim}}}
\newcommand{\Lim}{\mathrm{Lim}}

\newtheorem{theorem}{Theorem}[section]
\newtheorem{thm}[theorem]{Theorem}

\newtheorem{lem}[theorem]{Lemma}

\newtheorem{rem}[theorem]{Remark}
\newtheorem{ex}[theorem]{Example}

\begin{document}
\title{Corona limits of tilings : Periodic case}
\date{}
\author{Shigeki Akiyama, Jonathan Caalim, Katsunobu Imai,\\ and Hajime Kaneko}
\maketitle

\begin{abstract}
We study the limit shape of successive coronas of a tiling, which  
models the growth of crystals. We define basic terminologies and 
discuss the existence and uniqueness of corona limits, and then 
prove that corona limits are completely characterized by directional speeds. 
As an application, we give another proof that the corona limit of a 
periodic tiling is a centrally symmetric convex polyhedron (see \cite{Zhuravlev01, MalShu2011}). 
\end{abstract}

\section{introduction}


Motivated by a simple crystal growth model, 
we study the growth of coronas of a tiling of $\R^l$. The {\bf corona} of a patch (a finite set of tiles) 
is defined as the set of tiles sharing a boundary with the patch. Starting from a tile as a $0$-th corona, we can recursively define a sequence of $k$-th coronas ($k=0,1,2,\dots$). In this article, we study the {\bf corona limit}, 
a limit shape of the $k$-th corona as $k$ tends to $\infty$.

After submission, we are informed that corona limits were extensively researched by a Russian group of mathematicians including V.~G.~Zhuravlev, A.~V.~Shutov and A.~V.~Maleev, whose research work we are unaware of. 
They mainly studied the adjacency graph structure of tilings and described the step-by-step growth of the boundary of the coronas.
Many of the related papers are published in crystallography journals and some are currently only available in Russian.

Comparing with the existing literature, we find that our approach is a little more axiomatic. Indeed, we set up axioms for a general tiling and an adjacency relation defined on it that we can deal with. 
Then we introduce a key concept, {\bf directional speeds}, and show that the corona limit, if it exists, is completely described by these speeds as a star shape (Theorem \ref{Exist}). 
Under this setting, the convexity of a corona limit is deduced naturally from a certain uniformity of directional speeds (Theorem \ref{Uniform}).
We hope that this description gives a new insight into this study of corona limits and justify the raison d'\^{e}tre of this paper. 

In the case of periodic tilings, we reproduce
the result of \cite{Zhuravlev01, MalShu2011} that
the corona limit is a centrally symmetric convex polygon (Theorem \ref{periodic-coronalimit}). The shape of the corona limit depends on the shape of the $(3M-2)$-th corona where $M$ is the number of translationally inequivalent tiles. This fact leads to an effective algorithm to compute the corona limit of a periodic tiling. 


We examine our algorithm and compute the corona limits of all $1$-uniform and $2$-uniform tilings classified in \cite[Chapter 2.1-2.2]{Gruenbaum-Shephard:87} for both the point adjacency and edge adjacency. The corona limits are centrally symmetric convex polygons having 4, 6, 8, 10, 12 or 16 vertices.
We observe sensitivity to the adjacency, that is, the corona limit with respect to point adjacency may have a different shape from that of the limit obtained under the edge adjacency.  We also found a great variety of polygons: we are surprised that a decagon and a hexadecagon emerged from a periodic tiling, see Section \ref{Table}.
We do not know whether the number of vertices in a corona limit is bounded for general periodic tilings.

It is of great interest to generalize these results to non-periodic uniformly repetitive tilings\footnote{
Uniformity is necessary for the existence of corona limit. In fact, we give an example of a repetitive tiling witout a corona limit in Section \ref{NonUnif}.}. 
In a previous paper \cite{Akiyama-Imai:16}, the first and third authors showed that the corona limit of the Penrose tiling exists and it is
a regular decagon. The proof depends on a special property of the Penrose tiling. 
Shutov and Maleev \cite{ShuMal2017} gave another proof 
for this result. 
Growth shapes correspond to the boundaries of corona limits 
(for instance, see \cite{ShuMal2014}). 
The corona limit of a two dimensional Rauzy tiling, 
which is conjectured to be an octagon, was partially determined by 
Zhuravlev, Maleev \cite{ZhuMal2007} and 
Maleev, Shutov, Zhuravlev \cite{MalShuZhu2010}. 
The results in this direction are more or less example-driven and we are yet to reach comprehensive understanding.

\section{Basic definitions}
We study tilings in the Euclidean space $\R^l$. 
For a subset $X$ of $\R^l$, we denote its interior (resp. boundary) with respect to the Euclidean topology by 
$\mathrm{Inn}(X)$ (resp. $\partial(X)$). 
Let $\overline{X}$ be the closure of $X$. 
Moreover, let $\|\cdot\|$ denote the Euclidean norm on $\R^l$. 
A {\it tile} 
is a nonempty compact set 
which is the closure of its interior.  A tiling $\T$
is a covering of $\R^l$ by (countably many) tiles without interior overlaps. 
More precisely, 
a {\it tiling} $\T$ is a collection of tiles
$$
\{T_i \ |\ i=1,2,\dots\}
$$
having the properties
$$
\R^l=\bigcup_{i=1}^{\infty}  T_i
$$
and
$$
\mathrm{Inn} (T_i) \cap \mathrm{Inn}(T_j) \neq \emptyset
$$
implies $i=j$.  Denote by $B(x,r)$ the open ball  centered at $x$ of radius $r$. 
A tiling is {\it uniformly locally finite} if   
for any $r\in \R_{>0}$ there exists a positive 
integer $M_0$ such that for any $x\in \R^l$,  $B(x,r)$
is covered by at most $M_0$ tiles.
If $\T$ satisfies the axiom
\begin{itemize}
\item[(N)] There exist two positive real constants $s$ and $S$ such that for any tile $T$ in $\T$, there exists a ball of radius $s$ lying within $T$ and there exists a ball of radius $S$ which contains $T$.
\end{itemize}
in Gr\"unbaum-Shephard \cite{Gruenbaum-Shephard:87}, 
then $\T$ is uniformly locally finite. 
Indeed, let $B(x,r)$ be a ball such that there exist $M_0$ tiles $T_i$ with $T_i\cap B(x,r)\ne \emptyset$. 
Since every tile $T_i$ contains a ball of radius $s$, $B(x,r+2S)$ must contain $M_0$ disjoint balls of radius $s$. 
Comparing the volumes, we obtain
$$
\left(\frac {r+2S}s\right)^l\ge M_0.
$$
When 
$\T$ consists of finitely many tiles up to translation, the axiom (N) is obviously satisfied.

A subset $X$ of $\R^l$ is {\it relatively dense} if there exists a real constant $R>0$ such that
for any $x\in \R^l$, $B(x,R)\cap X\neq \emptyset$, and it is {\it uniformly discrete}
if there exists a real constant $r>0$ such that for any $x\in \R^l$, the set 
$B(x,r)\cap X$ is empty or a singleton. A {\it Delone set} is a subset of $\R^l$
which is both relatively dense and uniformly discrete. Delone sets model
atomic configuration of real materials. 
Choosing a suitable inner point in each tile of a tiling satisfying the axiom (N), we obtain a Delone set.
Conversely taking Vorono\"{i} partition (each point of $\R^l$ belongs to cells containing its closest points of $X$) of a Delone set $X$, we obtain a tiling satisfying the axiom (N). 
We say a tiling is of Delone type, or is a Delone tiling, if it satisfies the axiom (N). Throughout this paper, we assume that tilings are of Delone type.

It is well known that the family of nonempty compact sets of $\R^l$
forms a complete metric space by the Hausdorff metric:
$$
d_H(A,B)= \inf \{\eps \ge 0 \ |\ B \subset A[\eps] \ \mathrm{ and } \
A \subset B[\eps] \} 
$$
for compact subsets $A$ and $B$ of $\R^l$ with $K[\eps]=\{x\in \R^l\ |\
\text{there exists } k\in K \text{ that } \|x-k\| \le \eps \}$ where $K\subset \R^l$.

A patch $P$ of $\T$ is a nonempty finite subset of $\T$, i.e., 
$P=\{ T_i\ |\ i\in I\}$ with a nonempty finite subset $I$ of $\N$.
The support of $P$ is defined by $\supp(P)=\bigcup_{i\in I} T_i$.
An {\it adjacency} is a reflexive symmetric binary relation $\sim$ defined on $\T$
satisfying the following three conditions: 
\begin{enumerate} 
\item For any tiles $T_i$ and $T_j$, there exist
$n_1,\dots, n_{k}$ such that $i=n_1$, $j=n_k$ and $T_{n_{h}}\sim
T_{n_{h+1}}$ for $h=1,\dots, k-1$.
\item There exists a positive integer $M_1$ such that if $T_i\cap T_j\neq \emptyset$ 
then $k$ in a) can be chosen to be not greater than $M_1$.
\item There exists a positive real $R$ such that $T_i\sim T_j$ implies  
$d_H(T_i,T_j)\le R$. 
\end{enumerate}
 We say two tiles $T_i$ and $T_j$ in $\T$ are {\it adjacent} if 
$T_i \sim T_j$.  The relation $\sim$ induces an adjacency graph 
whose vertices are the elements of $\T$.

\begin{rem}\label{adjacency}
There are many ways to define adjacency.  
The relation $T_i\cap T_j\neq \emptyset$ defines the so-called {\it point adjacency}. 
If each $T_i$ is a 
manifold, then $\dim (T_i\cap T_j)\ge k$ defines {\it $k$-dimensional 
adjacency}. 
The case $k=1$ for $l=2$
is the edge adjacency used in \cite{Akiyama-Imai:16}.
\end{rem}

A patch $P$ 
is {\it connected} if
its induced subgraph is connected.  
Note that 
$\supp(P)$  is not necessarily topologically connected.
Given a patch $P$ and $n\in \N$, the $n$-th corona of $P$
is inductively defined as follows. The $0$-th corona is $P^{(0)}=P$ and the $(n+1)$-th corona
$P^{(n+1)}$ is the patch of $\T$
whose tiles are adjacent to some tile of $P^{(n)}$, i.e., 
$$
P^{(n+1)}=\{ T_i\in \T\ |\ \exists T_j\in P^{(n)} \quad T_i\sim T_j \}.
$$
In the following lemma, we consider the $n$-th corona of a patch $P=\{T_i\}$ consisting of a single tile $T_i$.

\begin{lem}\label{Linear}
Let $\T$ be a Delone tiling and $\sim$ an adjacency on $\T$. Then
there exist a real constant $c>0$ and a positive integer $n_0$ depending only on  
$\T$ and $\sim$ such that for all $i \in \N$ and $x \in T_i$, we have

\begin{equation}\label{eqnlin}
    B(x,n/c) \subset \supp(\{T_i\}^{(n)}) \subset B(x,cn)
\end{equation}
for any integer $n$ with $n\geq n_0$.
\end{lem}

\proof
For simplicity, put $T_k^{(n)}:=\supp(\{T_k\}^{(n)})$. Moreover, let $T:=T_i$ and $P=\{T\}$. 
The right-hand side of (\ref{eqnlin}) follows from 
$T^{(n)}\subset B(x,n(R+2S))$ for all $2\le n\in \N$ by axiom (N) and 
condition c) of $\sim$. 
On the other hand, 
by compactness, there exists a covering of $\partial(T^{(n)})[1/2]$ 
by a finite set of balls $B(x_k,1)$ with $x_k\in T^{(n)}$: 
$$\partial(T^{(n)})[1/2]\subset \bigcup_{k=1}^{\tau} B(x_k,1).$$ 
By uniform local finiteness of $\T$, each $B(x_k,1)$ is covered by 
at most $M_0$ tiles of $\T$. 
For any $1\leq k\leq \tau$, let $T_k\in P^{(n)}$ such that  $x_k\in T_k$. 
Then condition b) of $\sim$ implies that $B(x_k,1)\subset T_k^{(M_0M_1)}\subset T^{(n+M_0M_1)}$, and so 
$$
T^{(n)}[1/2] \subset T^{(n+M_0M_1)}.
$$
Let $x\in T$ and set $n=b_0 M_0M_1+b_1$ for integers $b_0,b_1$ with $1\leq b_1\leq M_0 M_1$. 
Using this inclusion relation successively, we obtain 
$$
B(x,cn)\supset
T^{(n)}\supset T^{(b_1)}[b_0/2]\supset B(x,b_0/2)\supset B(x,n/c)
$$
with $c=\max \{4M_0M_1, R+2S \}$ for $n\ge \max \{2, M_0M_1+1\}$. 
\qed

As a consequence, for any patch $P$ and point $x\in\supp(P)$, there is
a  constant $c$ such that
$$
B(x,n/c) \subset \supp(P^{(n)}) \subset B(x,cn),
$$
but clearly $c$ depends not only on $\T$ and $\sim$ but also on $P$. 
This suggests the study of the limit of the sequence of compact sets:
$$
\frac 1n \supp(P^{(n)})\quad n=1,2,\dots.
$$
Let us recall some basic facts on the convergence of compact sets. 
For nonempty compact sets $X_k\ (k=1,2,\dots)$, 
the {\it topological upper limit} $\LimU_{k} X_k$
is the set of all points $x$ such that every neighborhood of $x$ intersects infinitely many $X_k$'s,
and the {\it topological lower limit} $\LimL_k X_k$ is the set of all points $x$ such that
 every neighborhood of $x$ intersects $X_k$ for all but finite $k$. 
 Both the topological upper limit and topological lower limit are
 closed and we may rewrite the
 topological upper limit  $\LimU_k X_k=
 \bigcap_n \overline{\bigcup_{k\ge n} X_k}$. It is not empty 
if there exists a compact set $K$
that $X_k\subset K$.
If  $\LimL_k X_k=\LimU_k X_k$, the common value is called the {\it topological  limit} $\Lim_k X_k$.  
If the topological limit $\Lim_k X_k$ exists, then $X_k$ converges to $\Lim_k X_k$ 
by the Hausdorff metric. 
Conversely, if there exists a compact set $K$ such that $X_k\subset K$ 
and $X_k\rightarrow X$ by the Hausdorff metric, then
$X=\Lim_k X_k$  (see \cite{Kechris:95, Akiyama-Brunotte-Pethoe-Thuswaldner:08}).

Since there exists a positive constant $c$ such that $\supp(P^{(n)})/n \subset B(0, c)$, the
Hausdorff limit and the topological limit
coincide in our framework, when one of them exists.
We also see that $\LimL_k \supp(P^{(n)})/n \neq \emptyset$ from
$\supp(P^{(n)})/n \supset B(0, 1/c)$. 
We therefore define the {\it corona limit}:
$$
\lim_{n\rightarrow \infty} \frac 1n \supp(P^{(n)})
$$
by the Hausdorff metric, when the limit exists. 
Further we show that
if the corona limit exists, then it does not depend on the initial patch $P$. 

\begin{lem}
\label{Unique}
For two patches $P$ and $Q$ of a Delone tiling $\T$, we have
$$
\lim_{n\rightarrow \infty} \frac 1n \supp(P^{(n)})
=
\lim_{n\rightarrow \infty} \frac 1n \supp(Q^{(n)})
$$
if one of the limits exists.
\end{lem}
\proof
Within this proof,  we omit writing $\supp$.
By Lemma \ref{Linear}, there exists $M_2\in \N$ such that $Q\subset P^{(M_2)}$
and $P\subset Q^{(M_2)}$. Assume that 
$\lim_{n\rightarrow \infty} \frac 1n P^{(n)}$ exists. 
Then $$
\LimU_n  \frac 1n Q^{(n)} \subset \LimU_n \frac {n+M_2}n \frac 1{n+M_2}P^{(n+M_2)} 
=\LimU_n \frac 1n P^{(n)}=\lim_n \frac 1n P^{(n)}
$$
and
$$
\LimL_n \frac 1n Q^{(n)} =\LimL_n \frac 1{n+M_2} Q^{(n+M_2)} \supset
\LimL_n \frac n{n+M_2} \frac 1n P^{(n)}=\lim_n \frac 1n P^{(n)}.
$$
\qed

\begin{rem}
Shutov, Maleev, and Zhuravlev \cite{ShuMalZhu2009} proposed similar axioms on adjacency. Comparing with their axioms, we do not assume by our axioms that a transformation group acts on a tiling. 
Lemma \ref{Unique} is not new, see Theorem 1 in  \cite{ShuMal2014}). 
\end{rem}

Note that the origin has no particular geometric 
meaning in the definition of corona limits.
Indeed, $\lim_n \frac 1n P^{(n)}= \lim_n \frac 1n (P^{(n)}-x)$ for any $x\in \R^l$, 
one may think that any point in $\R^l$ is 
a center of growth of successive coronas.  

\section{Directional speed}
Let $\T=\{T_i\ |\ i=1,2,\ldots\}$ be a tiling. Let $\S=\{S_i\ |\ i=1,2,\ldots\}$ be a partition of 
$\R^l$. We say that $\S$ is supported by $\T$ if $\overline{S_i}=T_i$ for any $i\geq 1$. 
Note that if $\T$ is given, then there exists a partition $\S$ supported by $\T$ by the axiom of choice. 
By abuse of notation, we call  $S_i$ a tile in $\S$. We say that 
$P=\{S_{i_1},\ldots,S_{i_h}\}\subset \S$ is a connected patch in $\S$ if the corresponding patch 
$\{T_{i_1},\ldots,T_{i_h}\}$ is connected. Moreover, the support of $P$ is defined by 
$\supp(P)=\cup_{j=1}^h S_{i_j}$. \par
For $x,y\in \R^l$, we define a nonnegative integer $n(x,y)$ such that $1+n(x,y)$ is 
the minimum cardinality of connected patches in $\S$ whose support contains both
$x$ and $y$. 
Note that $n(x,y)=0$ if and only if $x,y$ belong to the same tile. 
Since $\S$ is a partition of $\R^l$, 
we have the inequalities:
\begin{equation}
\label{Triangle}
n(x,z)\le n(x,y)+n(y,z), \ 
|n(x,y)-n(x,z)|\leq n(y,z).
\end{equation}
We use the quantity $n(x,y)$ with a suitable partition $\S$ of $\R^l$ in order to investigate 
the asymptotic shape of $P^{(n)}$. 
The properties in Sections 3 and 4 do not depend on the choice of $\S$. \par
For a nonzero $\v\in \R^l$, we define the following four quantities to 
measure the speed of growth
along the direction $\v$:
\begin{eqnarray*}
\overline{d}_1(\v)&=& \limsup_{t\rightarrow \infty} \sup_{x\in \R^l} 
\frac {\| t\v\|}{n(x,x+t\v)}, \\
\overline{d}_2(\v)&=& \sup_{x\in \R^l} \limsup_{t\rightarrow \infty} 
\frac {\| t\v\|}{n(x,x+t\v)},\\
\underline{d}_2(\v)&=& \inf_{x\in \R^l} \liminf_{t\rightarrow \infty} 
\frac {\| t\v\|}{n(x,x+t\v)},\\
\underline{d}_1(\v)&=& \liminf_{t\rightarrow \infty} \inf_{x\in \R^l} 
\frac {\| t\v\|}{n(x,x+t\v)}. 
\end{eqnarray*}

Let $d$ be one of $\overline{d}_i,\underline{d}_i$ (i=1,2).
By definition we have 
$d(r\v)=d(\v)$ for all $r>0$ and $d(\v)$ is determined by the values on the $l-1$
unit sphere.

\begin{lem}
\label{Cont}
For a Delone tiling $\T$, the four quantities $\overline{d}_i, \underline{d}_i$ for $i=1,2$ are 
continuous on $\v\in \R^l$, taking finite and positive values.
\end{lem}

\proof
We prove the finiteness, positivity, and continuity as a function of the variable $\v$ with $\|\v\|=1$. Take a partition $\S=\{S_i\ |\ i=1,2,\ldots\}$ supported by $\T$. 
Choose $i$ such that $x\in S_i$. For a sufficiently large $t$, we have
$$
x+t \v \in T_i^{(n)}\setminus T_i^{(n-1)}
$$
with $n=n(x,x+t\v)\ge 2$. 
Using Lemma \ref{Linear},
there are positive constants $C_0,C_1$ depending on $\T$, $\sim$ such that
for all $x\in \R^l$, $\|\v\|=1$ and sufficiently large $t$, 
\begin{equation}
\label{Linear2}
C_0 t\le n(x,x+t\v)\le C_1 t.
\end{equation}
This implies that the four quantities in Lemma \ref{Cont} are positive and finite.
Therefore for all $x\in \R^l$, $\v\in \R^l$ there exists $t_0>0$ that if $t>t_0$, then  
$$
C_0 t \|\v\|\le
n(x,x+t\v)\le C_1 t \|\v\|.
$$
Let $t$ be sufficiently large. Using (\ref{Triangle}), for any $\v,\v'$ with $\|\v\|=\|\v'\|=1$, we have 
\begin{eqnarray*}
&& \left|\frac {\|t\v'\|}{n(x,x+t\v')}-\frac{\|t\v\|}{n(x,x+t\v)}\right| \\
&\le&
\frac{1}{C_0^2 t}\big|n(x,x+t\v)-n(x,x+t\v')\big|\\
&\le& 
\frac{1}{C_0^2 t}\bigl(n(x+t\v,x+t \v+t(\v'-\v))
     \bigr)\\
     &\le&\frac{C_1}{C_0^2}\|\v-\v'\|,
\end{eqnarray*}
which gives the required continuity.
\qed
\medskip

\begin{rem}
\label{LinearUpper}
Lemma \ref{Linear} immediately implies the right inequality of (\ref{Linear2}) for all $t\ge 1$, i.e., 
there exists a positive constant $C_1$ that
for any $x\in \R^l, \|\v\|=1$ and $t\ge 1$, we have $n(x,x+t\v)\le C_1 t$.
\end{rem}

By Lemma \ref{Cont},
$d(\v)$ is uniformly continuous in $\R^l\setminus \{0\}$. 
We prove the following inequalities.

\begin{lem}For all $0\neq \v \in \R^l$, we have
\label{Basic}
$$
\overline{d}_1(\v) \ge 
\overline{d}_2(\v) \ge 
\underline{d}_2(\v) \ge 
\underline{d}_1(\v) . 
$$
\end{lem}

\proof
The middle inequality is obvious.  Put $g(x,t)=t\| \v\|/n(x,x+t \v)$.
Assume that
$\overline{d}_1(\v) <\overline{d}_2(\v)$. 
Then for any $\eps>0$, 
there exists $t_0 \in \R$ such that for any $t\ge t_0$ and any $x\in \R^l$, 
$g(x,t)<\overline{d}_1 (\v)+\eps$.
On the other hand, for any $\eps>0$
there exist $x_0\in \R^l$ and infinitely many $t$'s such that 
$g(x_0,t)> \overline{d}_2(\v)-\eps$. 
Choosing $\eps=(\overline{d}_2(\v)-\overline{d}_1(\v))/3$,
we have an impossibility
$$
g(x_0,t)\ge \overline{d}_2(\v)-\eps > \overline{d}_1(\v)+\eps\ge g(x_0,t)
$$
for some $x_0\in \R^l$ and some $t\ge t_0$. 
The proof of the last inequality is similar\footnote{
F.~Nakano suggested us an easier proof of Lemma \ref{Basic}, which works for any real valued function $g(x,t)$ in two variables.}.
\qed

\begin{ex}
\label{12Tiling}
Consider the tiling of $\R$ by $T_1=[0,1]$ and $T_2=[0,2]$ where  
the positive part of $\R$, i.e. $[0,\infty)$, is tiled by translates of $T_2$
and the negative part of $\R$, i.e. $(-\infty, 0]$, is tiled by translates of $T_1$ in an obvious manner. 
Then we see
$1=\underline{d}_1(1)<\underline{d}_2(1)=\overline{d}_2(1)=\overline{d}_1(1)=2$.
For the negative direction, we have
$1=\underline{d}_1(-1)=\underline{d}_2(-1)=\overline{d}_2(-1)<\overline{d}_1(-1)=2$.
\end{ex}

\begin{ex}
Using the same notation as in Example \ref{12Tiling}, we tile the negative part of $\R$ by translates of $T_1$. 
Decompose $\R_{\ge 0}$ into $[a_0,a_1]\cup [a_1,a_2] \cup \dots$ by an 
increasing even integer
sequence $(a_n)$ with $a_0=0$. We tile 
$[a_{2n},a_{2n+1}]$ by translates of $T_2$ and $[a_{2n-1},a_{2n}]$ by translates of $T_1$.
Choosing the sequence $(a_n)$ of rapid growth, we can
construct an example such that $\underline{d}_2(1)<\overline{d}_2(1)$, e.g., 
taking $a_n=2^{2^n}-2$, we have $\underline{d}_2(1)=1 < \overline{d}_2(1)=2$.
\end{ex}

\begin{figure}[h]
\begin{center}
\includegraphics[scale=0.20]{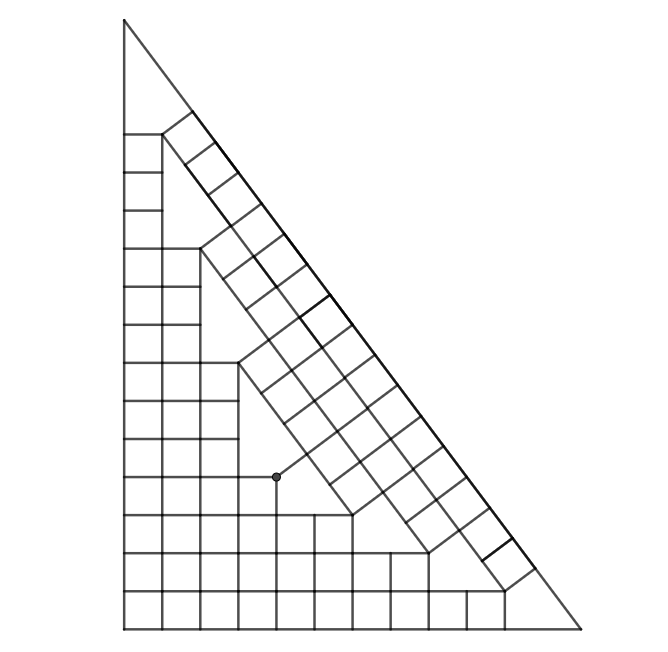}
\caption{The tiling of $\R^2$ described in  Example~\ref{345}.}
\label{caalim-tiling}
\end{center}
\end{figure}

\begin{ex}
\label{345}
Consider a right triangle $\Delta$ whose edges have length 3,4 and 5,
which has the inscribed circle of radius 1. Assume that its incenter is located at the origin and the edge of length 3 is parallel to the $x$-axis, and the edge
of length 4 is parallel to the $y$-axis.
Then $(n+1)\Delta-n \Delta$ is subdivided using unit squares placed along the edges of $n\Delta$ and two quadrangles appear at two non-right interior angles of $\Delta$. This gives a tiling of $\R^2$ by 4 translationally inequivalent tiles (Figure~\ref{caalim-tiling}).
By construction, we immediately see that the corona limit by point adjacency is $\Delta$ itself. 
For $\v=(4,3)$, we see $$5/4=\overline{d}_1(\v)>\overline{d}_2(\v)=\underline{d}_2(\v)=\underline{d}_1(\v)=1,
$$ and for $\v=(-1,-1)$,
$$\sqrt{2}=\overline{d}_1(\v)=\overline{d}_2(\v)=\underline{d}_2(\v)>\underline{d}_1(\v)=1.$$
\end{ex}

These examples suggest that a strict inequality in Lemma \ref{Basic} came from a certain lack of uniformity of the corresponding tiling. See Section \ref{NonUnif} and 
the third problem in Section \ref{Final}.

\begin{lem}[Directional Speed]
\label{DSpeed}
Let $\v\in \R^l$ with $\|\v\|=1$. The following conditions are equivalent:
\begin{enumerate}
\item for any $\eps>0$ there exists a relatively dense set $S=S(\v,\eps)
\subset \R^l$ 
such that for any $x\in S$ we have
$$
\liminf_{t\rightarrow \infty} \frac {t \|\v\|}{n(x,x+t\v)}>\overline{d}_2(\v)-\eps,
$$
\item $\frac {t \|\v\|}{n(x,x+t\v)}$ converges to a value 
independent of $x\in \R^l$ as $t$ tends to infinity.
\item $\overline{d}_2(\v)=\underline{d}_2(\v)$.
\end{enumerate}
If one of the conditions holds, we say 
that $\T$ admits a directional speed in the direction $\v$, or $d(\v)$ exists, where 
$d(\v)=\overline{d}_2(\v)$.
\end{lem}

\proof
We show that a) implies c). Let $\eps>0$. There exists 
a constant $R>0$ such that for any $y\in \R^l$ we find $x\in S$ with $\| x-y\|\le R$. 
By assumption, we find $t_0=t_0(x)>0$ such that if $t\ge t_0$, then
\begin{equation}
\label{Below1}
\frac {t \|\v\|}{n(x,x+t\v)} > \overline{d}_2(\v)-\eps.
\end{equation}
Remark \ref{LinearUpper} shows $|n(y,y+t\v)-n(x,x+t\v)|< 2 C_1 R$. 
Thus we also have
$$
\frac {t \|\v\|}{n(y,y+t\v)} > \overline{d}_2(\v)-2\eps.
$$
 From the definition of $\underline{d}_2$, for any $y\in \R^l$, there exists $t_1=t_1(y)>0$ such that 
\begin{equation}
\label{Above1}
\frac {t \|\v\|}{n(y,y+t\v)} < \underline{d}_2(\v)+\eps.
\end{equation}
for $t\ge t_1$. Combining (\ref{Below1}) and (\ref{Above1}), we see c). The equivalence
of b) and c) is clear, and b) implies a) by taking $S=\R^l$.
\qed

\begin{thm}
\label{Exist}
Let $\T$ be a Delone tiling and $\sim$ an adjacency
on $\T$. Then the corona limit $K=\lim_n \supp(P^{(n)})/n$ exists 
if and only if $\T$ admits
a directional speed in every direction 
$\v\in \R^l$ with $\|\v\|=1$.
Further we have
\begin{equation}
\label{Star}
K=\left\{ t d(\v)\v \ |\  t\in [0,1], \|\v\|=1 \right\}.
\end{equation}
\end{thm}

\proof
By Lemma \ref{Unique} and the remark after it, we may assume that 
$P=\{T\}$ and $0\in \mathrm{Inn} (T)$. 
Let $\v\in \R^l$ with $\| \v\|=1$. By the compactness of $T^{(m)}$, we can define $\sigma(m)=\max \{t\ |\ t\v \in T^{(m)} \}$ . 
From Lemma \ref{Linear}, we have $m/c<\sigma(m)<cm$ for some positive constant $c\in \R$.

Assume that $\lim_m T^{(m)}/m$ exists. 
Then we see 
$
\lim_{m\rightarrow \infty} \frac{\sigma(m)}m
$
exists and is equal to $\lim_m \frac{\sigma(m)}{n(0,\sigma(m)\v)}$.
By a similar argument as in the proof of Lemma \ref{DSpeed},
we get for all $x\in \R^l$ that 
$\frac {t \|\v\|}{n(x,x+t\v)}$ converges to $\lim_m \frac{\sigma(m)}{n(0,\sigma(m)\v)}$.

On the other hand, suppose that $d(\v)$ exists for every $\v\in \R^l$ with $\| \v\|=1$. 
Then for sufficiently large $m\in \N$, we have
$$
\sigma(m)<(d(\v)+\eps)n(0,\sigma(m)\v)=(d(\v)+\eps)m,
$$
and so 
\begin{equation}
\label{aaa}
\{ t\v \ |\ t\ge 0\} \cap \LimU_m \frac{1}{m} T^{(m)}
\subset
\{ t d(\v)\v\ |\ t\in[0,1]\}.
\end{equation}
Moreover, let $0<t<1$. For any $\eps>0$, we get that 
$n(0,m t d(\v) \v)<(t+\eps)m$ for any sufficiently large $m$, and so 
\begin{equation}
\label{bbb}
t d(\v) \v \in \{ t\v \ |\ t\ge 0\} \cap \LimL_m \frac{1}{m} T^{(m)}.
\end{equation}
Since (\ref{aaa}) and (\ref{bbb}) hold for any $\v$ and $0<t<1$, we obtain that $\lim_m T^{(m)}/m$ exists 
and (\ref{Star}) is satisfied. 
\qed

By Lemma \ref{Cont}, we know that if the directional speed $d(\v)$ 
exists then it is finite, positive and continuous with respect to the variable $\v$.
This implies $\partial(K)=\left\{ d(\v)\v \ |\ \|\v\|=1 \right\}$.
Conversely if $K$ is a compact star convex set written as
$$
K=\left\{ t g(\v)\v \ |\  t\in [0,1], \|\v\|=1 \right\}
$$
with a continuous positive function $g$ defined on the $l$-dimensional unit sphere, 
then it is realized as the corona limit of a 
Delone tiling where the adjacency is given by point adjacency. The construction of such tiling is
similar to that of  Example \ref{345}, i.e., we dissect 
each annulus $(n+1)K-nK$ into 
small pieces by hyperplanes, chosen appropriately to satisfy axiom (N). 

Therefore the corona limit is characterized as a star shape centered at the origin, containing the origin in its interior having a continuous gauge function.

\section{Shapes of corona limits}

\begin{lem}[Uniform Directional Speed]
\label{UDSpeed}
Let $\v\in \R^l$ with $\|\v\|=1$. The following conditions are equivalent.
\begin{enumerate}
\item $\frac {t \|\v\|}{n(x,x+t\v)}$ converges uniformly with respect to $x\in \R^l$ to a value 
independent of $x\in \R^l$ as $t$ tends to infinity.
\item $\overline{d}_1(\v)=\underline{d}_1(\v)(=d(\v))$.
\end{enumerate}
If one of the conditions holds, we say 
that $\T$ admits a uniform directional speed in direction $\v$.
\end{lem}

\proof
We see that b) is equivalent to  the fact that 
for any $\eps>0$, there exists $t_0$ such that for any $t\ge t_0$,  
$$
d(\v)-\eps < \frac {t \|\v\|}{n(x,x+t\v)} <d(\v)+\eps
$$
holds for any $x\in \R^l$, i.e, the uniformity of the convergence. The converse is also clear.
\qed

\begin{thm}
\label{Uniform}
Let $\T$ be a Delone tiling and $\sim$ an adjacency on $\T$. If 
$\T$ admits
a uniform directional speed in every direction, then the corona limit
is convex and symmetric with respect to the origin.
\end{thm}

\proof
By Lemma \ref{UDSpeed}, for any $\eps>0$ there exists $t_0$ 
such that if $t\ge t_0$ then for any $x\in \R^l$, 
$$
\left| \frac {t\|\v\|}{n(x,x+t\v)}-d(\v)\right| <\eps.
$$
Let $y\in\R^l$ and put $x=y-t\v$. We have
$$
\left| \frac {t\|\v\|}{n(y,y-t\v)}-d(\v)\right| <\eps,
$$
for any $t\ge t_0$ because $n(a,b)=n(b,a)$. Thus the corona limit is symmetric with respect to
$\v \leftrightarrow -\v$, and we have
\begin{equation}\label{dvseta}
d(\v)= \lim_{t\rightarrow \infty} \frac{\|t\v\|}{\eta(t\v)}
\end{equation}
with 
\begin{equation}\label{defeta}
\eta(\v) = \min \{n(x,y)\ |\ x,y\in \R^l, \ x-y=\v\}.
\end{equation} 
We have $\eta(\v)=\eta(-\v)$ and $d(\v)=d(-\v)$. 
Therefore the corona limit is symmetric with respect to the origin. 

Since $K$ is a star shape by Theorem \ref{Exist},   
it is enough to show that for all $\v_1, \v_2\in \partial(K)$ and $u\in (0,1)$ with $u\v_1+(1-u)\v_2\ne 0$, we have 
$$d(u\v_1+(1-u)\v_2)\ge \|u\v_1+(1-u)\v_2\|$$
to prove the convexity of $K$.
Let $\v=u\v_1+(1-u)\v_2$. For $x\in \R^l$, put 
$A(t)=n(x,x+tu \v_1)$ and $B(t)=n(x+tu \v_1,x+tu \v_1+t(1-u)\v_2)$. Hereafter 
we write $f(t)=o(g(t))$ to mean that $\lim_{t\rightarrow \infty}\|f(t)\|/g(t)=0$, which is
the Landau symbol on the vector valued function $f(t)$.
Observe that the following three equations hold:
\begin{eqnarray*}
&&n(x,x+t\v)d(\v)\v = t\|\v\| \v + o(t),
\\
&&A(t)d(\v_1)\v_1 = tu\|\v_1\| \v_1 + o(t),
\\
&&B(t)d(\v_2)\v_2  = t(1-u)\|\v_2\| \v_2 + o(t).
\end{eqnarray*}
Recall for $i=1,2$ that $d(\v_i)=\|\v_i\|$ because $\v_i\in \partial(K)$. 
Setting $\v'=(1/\|\v\|) \v$, we obtain
$$
\frac {n(x,x+t\v)}{A(t)+B(t)}d(\v')\v' = \left(\frac{A(t)}{A(t)+B(t)} \v_1 + \frac {B(t)}{A(t)+B(t)} \v_2 \right) +o(1).
$$
From $A(t)+B(t)\ge n(x,x+t\v)$, letting $t$ tend to infinity,
we have
$$
\kappa d(\v')\v' = u\v_1+(1-u)\v_2=\v
$$
with $0<\kappa\le 1$. This implies
$d(\v)\ge \|\v\|$.
\qed

\section{Periodic tilings}\label{Period}
In applications, we often consider the case where a tiling $\T$ consists of
finitely many tiles up to translations, or up to more general rigid motions. We denote by $G$ such a transformation 
group acting 
on $\T$ and by $g(T)\in \T$ the image of a tile  $T$ by $g\in G$. Shutov, Maleev, and Zhuravlev \cite{ShuMalZhu2009} proposed an assumption that 
for any $T_i, T_j\in \T$ and any $g\in G$, we have 
$T_i\sim T_j$ if and only if $g(T_i)\sim g(T_j)$.
Hereafter we consider the case where $G$ acts as translations. 
A translation of $\T$ by $x\in \R^l$ is defined by $\T+x=\{T_i+x\ |\ i\in \N\}$.
We say that $x\in \R^l$ is a period of $\T$ if $\T=\T+x$.  
Let $\L$ be the set of periods of $\T$.
A tiling is {\it lattice periodic} if there exist $l$ periods which are linearly
independent over $\R$, i.e., $\L$ forms a lattice in $\R^l$. Clearly there are only finitely many tiles up to translations in $\T$, and $G=\L$ acts on the tiling. 
In what follows, we assume that for any $p\in \L$, we have $T_i\sim T_j$ if and only if $T_i+p\sim T_j+p$. In other words, we only care about translational symmetry of the tiling.

We now need an additional technical assumption. Let $\S=\{S_i \ | \ i=1,2,\ldots\}$ be a partition of 
$\R^l$ supported by a lattice periodic tiling $\T$ with set $\L$ of periods. 
By the axiom of choice, we may assume that $\S$ satisfies the following: 
For any $p\in \L$ and $S_i\in \S$, we have $p+S_i\in \S$. 
Under this assumption, we get for any $p\in \L\setminus \{0\}$ that $\eta(p)>0$, where $\eta(\v)$ is defined by (\ref{defeta}). 
In fact, for any $S_i \in \S$ and $y\in S_i$, we see that 
$n(p+y,y)>0$ because $p+y\not\in S_i$ by $p+y\in p+S_i\in \S$. 

\begin{thm}\label{periodic-coronalimit}
Let $\T$ be a lattice periodic tiling. Then the corona limit of $\T$ is
a convex polyhedron which is symmetric with respect to the origin.
\end{thm}

\proof
Fix a fundamental domain $F$ of $\R^l/\L$ whose closure is 
compact. Take $\v\in\R^l$. The supremum of the function 
$
f(x):=\limsup_{t\rightarrow \infty} \frac{\| t \v\|}{n(x,x+t\v)}
$
is attained at a point $x_0\in F$ because $f(x)$ is continuous. 
Denote by $\lfloor {\bf x} \rfloor$ the unique element of
$\L$ such that $\{\mathbf{x}\}:=\mathbf{x}-\lfloor \mathbf{x} \rfloor\in F$. 
Put $a_m=n(x,x+m \v)$ for $m \in \N$. Then there exists a positive
constant $c$ such that
$$
a_{m+m'}\le a_m+a_{m'}+c.
$$
In fact, we have
\begin{eqnarray*}
n(x,x+(m+m')\v)& \le& n(x,x+m\v)+n(x+m\v,x+(m+m')\v)  \\
& = & n(x,x+m\v)+n(x+\{m\v\}, x+\{m\v\} + m'\v)\\
& \le& n(x,x+m\v)+n(x,x+m'\v)+c
\end{eqnarray*}
where $c=2 \max_{a,b\in 2F} n(a,b)$. 
Since $a_m+c$ is subadditive, a well-known principle (c.f. \cite[Theorem 4.9]{Walters})
implies
$$
\lim_m \frac{a_m}m=\lim_m \frac{a_m+c}m=\inf_m \frac {a_m+c}m=\inf_m \frac{a_m}m.
$$
Switching to positive real variable $t$ is plain, and we see 
that $\lim_{t\rightarrow \infty} \frac{t\|\v\|}{n(x,x+t\v)}$ exists
for each $\v\in \R^l\setminus \{0\}$. 
Theorem \ref{Exist} implies that the corona limit $K$ exists because $S=x_0+\L$ satisfies the condition a) of Lemma \ref{DSpeed}. 
Moreover the convergence is uniform with respect to $x\in \R^l$ because there exists a constant 
$c'>0$ such that, for any $x\in \overline{F}$ and $\v \in \R^l$, 
$$
|n(x,x+t\v)- n(0,t\v)|\le n(0,x)+n(t\v, x+t\v ) \le c'
$$
by the uniform local finiteness of $\T$. 
The inequality above also holds for any $x\in \R^l$ by the periodicity of $\T$. 
Theorem \ref{Uniform} implies that $K$ is compact, convex, and symmetric with respect to the origin.

In what follows, we show that $K$ is a polyhedron. 
We say that $x,y\in \R^l$ are $\L$-equivalent if $x-y\in \L$. 
Similarly, we define the $\L$-equivalence of two tiles. 
We say that $x\in \L$ is {\it decomposable} if there exists a decomposition 
$x=x_1+x_2$ with $x_1,x_2\in \L\setminus \{0\}$ satisfying 
$\eta(x)\ge \eta(x_1)+\eta(x_2).$ 
Let $M$ be the number of inequivalent tiles under $\L$-equivalence. 

We now verify that if $x\in\L$ satisfies $n(x,0)>3M-2$, then $x$ is decomposable. 
We claim that $\eta(x)>M$. 
We use this claim to ensure that a connected patch $Q$ defined later contains translationally equivalent tiles. 
Suppose on the contrary that $\eta(x)\leq M$. 
Then there exist $y_1,y_2\in \R^l$ satisfying $y_1-y_2=x$ and $\eta(x)=n(y_1,y_2)$. 
We may assume that $y_1$ satisfies the minimality:
$$
n(0,y_1)=\min\{n(0,y')\ | \ y'-y_1\in \L\}.
$$
Let $P$ be a connected patch in $\S$ satisfying $0,y_1 \in \supp(P)$ and 
$|P|=1+n(0,y_1)$, that is, 
$P$ is a connected patch in $\S$ with minimum cardinality satisfying $0,y_1 \in \supp(P)$. 
Recall that $\L$ is also a period of $\S$. Thus, the minimality on $y_1$ 
implies that any two tiles in $P$ are not $\L$-equivalent, and so $n(0,y_1)\leq M-1$. 
Since $y_1-y_2\in\L$, we see 
$$n(y_1-y_2,y_1)=n(0,y_2)\leq n(0,y_1)+n(y_1,y_2)\leq 2M-1,$$
and so $n(x,0)\leq 3M-2$. This shows the claim. 

There exists $y\in \R^l$ with $\eta(x)=n(y+x,y)$.
Let $Q$ be a connected patch in $\S$ satisfying $y+x,y \in \supp(Q)$ and 
$|Q|=1+n(y+x,y)(\ge M+2)$. 
Excluding the tile containing $y$, there are two translationally
equivalent tiles by $|Q|-1\geq M+1$. 
Pick two $\L$-equivalent points $z_1,z_2$ from the interiors of these two tiles. 
By the minimality of $Q$, by switching indexes if necessary, 
the adjacency graph induced by $\sim$ on $Q$ is a linear graph and 
the patch $Q$ is divided into three patches $A,B$ and $C$ such that 
$y+x,z_1 \in \supp(A)$, $z_1,z_2 \in \supp(B)$ and $z_2,y \in \supp(C)$. 
Note that $y\neq z_2$ and the linear structure and the minimality of $Q$
imply $x\neq z_1-z_2$. Using the minimality again, we see 
$$
\eta(x)=n(y+x,y)=n(y+x,z_1)+n(z_1,z_2)+n(z_2,y).
$$
Since $z_1,z_2$ are $\L$-equivalent, $C-(z_2-z_1)$ is a patch in $\T$
which connects $z_1$ and $y-(z_2-z_1)$ and $\T=\T-(z_2-z_1)$.
Thus, we obtain 
\begin{eqnarray*}
\eta(x)&=& n(y+x,z_1)+n(z_1,z_2)+n(z_1,y-(z_2-z_1))\\
&\ge& n(y+x,y-(z_2-z_1))+n(z_1,z_2)\\
&\ge& \eta(x-(z_1-z_2))+\eta(z_1-z_2).
\end{eqnarray*}
Hence, $x$ is decomposable because $x-(z_1-z_2)$ and $z_1-z_2$ are in $\L\setminus \{0\}$. 

Observe that 
\begin{eqnarray*}
\L'&:=&\{x\in \L\setminus\{0\} \ | \ n(x,0)\leq 3M-2\},\\
\F&:=&\{(1/\eta(x)) x\ | \ x\in \L'\}
\end{eqnarray*}
are finite sets because $\|x\|\leq (3M-2)R$ for any $x\in \L'$. 
Thus, the convex hull $K'$ of $\F$ is a polyhedron. 
Our goal is to show that $K=K'$. For the proof of $K'\subset K$, it suffices to check for any $x\in \L'$ that $d(x)/\|x\|\geq 1/\eta(x)$. There exists $y\in \R^l$ with $n(y+x,y)=\eta(x)$. Observing for any $m\geq 1$ that $n(y+mx,y)\leq m\eta(x)$ by the periodicity of $\T$, we get 
$$
\frac{d(x)}{\|x\|}=\lim_{m\to\infty} \frac{m}{n(y+mx,y)}\geq \frac1{\eta(x)}.
$$
For the proof of $K\subset K'$, it suffices to verify for any $\v\in \R^l$ with $\|\v\|=1$ that 
$d(\v)\v\in K'$. 
We consider the sequence $\lfloor m\v\rfloor$ ($m=1,2,\ldots$). For any sufficiently large $m$, we see $\lfloor m\v\rfloor$ is decomposable by $n(\lfloor m\v\rfloor,0)>3M-2$. Denote the corresponding decomposition of $\lfloor m\v\rfloor$ by $\lfloor m\v\rfloor=x_1+x_2$ with $x_1,x_2\in\L\setminus\{0\}$. 
If $n(x_1,0)>3M-2$ (resp. $n(x_2,0)>3M-2$), then we construct a successive decomposition of $x_1$ (resp. $x_2$). The decomposition terminates in  finite time because $\eta(y)\geq 1$ for any $y\in \L\setminus\{0\}$. 
Hence, we get a decomposition
$$
\lfloor m\v\rfloor=\sum_{x\in \L'} c(x,m) x
$$
with 
$$
\eta(\lfloor m\v\rfloor)\ge \sum_{x\in \L'} c(x,m) \eta(x),
$$
where $c(x,m)$ is a nonnegative integer for any $x\in \L'$. 
Putting 
\begin{eqnarray*}
t:=\frac {\sum_{y\in L'} c(y,m)\eta(y)}{\eta(\lfloor m\v\rfloor)}\in[0,1], \ 
u(x):=\frac{c(x,m)\eta(x)}{\sum_{y\in \L'} c(y,m) \eta(y)}\in [0,1]
\end{eqnarray*}
for $x\in L'$, we obtain 
$$
\frac{1}{\eta(\lfloor m\v\rfloor)}\lfloor m\v\rfloor=t \sum_{x\in L'} u(x)\cdot \frac{1}{\eta(x)} x\in K'
$$
by $\sum_{x\in L'} u(x)=1$. 
Note that $\{m\v\}$ ($m=1,2,\ldots$) is bounded and that (\ref{dvseta}) holds because the assumption of Theorem 
\ref{Uniform} is satisfied. Therefore, we deduce from Lemma \ref{Linear} and $\|\v\|=1$ that 
$$
\lim_{m\rightarrow\infty} \frac1{\eta(\lfloor m\v\rfloor)}\lfloor m\v\rfloor
=\lim_{m\rightarrow \infty} \frac{m}{\eta(m\v)}\v=d(\v)\v\in K'
$$ 
because $K'$ is closed. 
\qed
\medskip

For a compact convex set $X$, an element $x\in X$ is an {\it extreme} point
if $x=u a+(1-u)b$ for $a,b\in X$ and $u\in [0,1]$ implies $x=a$ or $x=b$.
Minkowski-Carath\'eodory
Theorem (\cite{Simon})\footnote{An infinite dimensional version is due to
Krein-Milman (cf. \cite{Simon,Conway}).} implies that $X$ is the 
convex hull of the set of extreme points, in particular, the set of
extreme points is non empty. 
The above proof shows that for any lattice periodic tiling, the corona limit
has finitely many extreme points, i.e., the vertices of $K$.
Moreover we have shown that 
some integer multiple of each extreme point becomes a period of $\L$. This implies that the set of extreme points lies in the $l$-dimensional vector space over $\Q$, that is,  
\cite[(iii) Theorem 3]{ShuMal2014}. This answers the second question in Section \ref{Final}.

For actual computation, we have a simple way to compute corona limits.  
Collect all connected linear patch consisting of $k(\le M+1)$ tiles
whose two tiles at the ends are $\L$-equivalent by a period $\v\in \L$. 
Then the convex hull of these $\v/k$'s is the desired corona limit. We call these $\v/k$'s {\bf velocities} in Section \ref{Table} and Appendix~\ref{sec:1-regular-plots} and \ref{sec:2-regular-plots}.

\section{Corona limits of $k$-uniform tilings}
\label{Table}

From the proof of Theorem~\ref{periodic-coronalimit}, 
we can devise an algorithm that computes the velocities and all the extreme points of a corona limit. In this section, 
the terms corona limit and edge corona limit refer to the limits with respect to the point adjacency and edge adjacency, respectively.
A tiling is $k$-uniform if it is tiled by regular polygons and has $k$ different vertex configurations 
under the action of its symmetry group, see
\cite[Chapter 2.1-2.2]{Gruenbaum-Shephard:87}.
The $1$-uniform tilings are often called Archimedean tilings and are individually named according to the vertex configurations.
For example, the Archimedean tiling $3^2.4.3.4$ 
means that in the tiling the tiles of shape triangle, triangle, square, triangle and square surround every vertex in this cyclic order.
Appendix~\ref{sec:1-regular-plots} and \ref{sec:2-regular-plots}  list the plots of the velocities of 1-uniform and 2-uniform tilings, respectively. Here, the edge length of each regular polygon is  
normalized to $1$.

Sensitivity to the adjacency is already observed in
the square tiling $4^4$: both corona limits are squares but the sizes and edge directions are different.
Moreover, corona limits of the same tiling may not be affine equivalent. For the tiling $3^2.4.3.4$ ((1-09) in Appendix~\ref{sec:1-regular-plots}), 
the corona limit has four extreme points:
$$
\pm \left(\frac{3+\sqrt{3}}4, \frac{1+\sqrt{3}}4\right),
\pm \left(\frac{-1-\sqrt{3}}4, \frac{3+\sqrt{3}}4\right)
$$
forming a square.
Meanwhile, the edge corona limit has eight extreme points as described below. It is not a regular octagon (see Section \ref{Final}), nevertheless, the edge lengths are all equal to $\sqrt{10(2+\sqrt{3})}/12$.

Most of the shapes are quadrilaterals or hexagons except the following 6 cases listed with their extreme points.   
\begin{itemize}
\item Octagon: edge corona limit of $3^2.4.3.4$ (1-09)\\
\begin{equation*}
\pm \left(\frac{3+\sqrt{3}}8, \frac{1+\sqrt{3}}8\right), \pm \left(\frac 16, \frac{2+\sqrt{3}}6\right), \pm \left(\frac{-1-\sqrt{3}}8, \frac{3+\sqrt{3}}8\right), \pm \left(\frac{-2-\sqrt{3}}6, \frac 16\right)
\end{equation*}
\item Octagon: corona limit of $3^3.4^2;3^2.4.3.4$ (2-15) \\
\begin{equation*}
\pm\left(\frac{3+\sqrt{3}}4, 0\right), 
\pm\left(\frac{3+\sqrt{3}}6, \frac{3+\sqrt{3}}6\right),
\pm\left(0,\frac{3+\sqrt{3}}4\right),
\pm\left(\frac{-3-\sqrt{3}}6, \frac{3+\sqrt{3}}6\right)
\end{equation*}
\item Octagon: corona limit of $3^3.4^2;3^2.4.3.4$ (2-16)\\
\begin{equation*}
\pm\left(\frac{5+2 \sqrt{3}}{8}, \frac{\sqrt{3}}{8}\right),
\pm\left(\frac{7+2 \sqrt{3}}{12} , \frac{4+3 \sqrt{3}}{12}  \right),
\end{equation*}
\begin{equation*}
\pm\left(\frac{-1-2 \sqrt{3}}{12}, \frac{8+3 \sqrt{3}}{12} \right),
\pm\left(\frac{-3-2 \sqrt{3}}{8}, \frac{4+\sqrt{3}}{8} \right)
\end{equation*}
\item Decagon: edge corona limit of $3^3.4^2;3^2.4.3.4$ (2-16) \\
\begin{equation*}
\pm\left(\frac{5+2 \sqrt{3}}{12} , \frac{\sqrt{3}}{12} \right),
\pm\left(\frac{3+\sqrt{3}}{8} , \frac{1+\sqrt{3}}{8} \right),
\pm\left( \frac{1}{6}, \frac{2+\sqrt{3}}{6} \right),
\end{equation*}
\begin{equation*}
\pm\left( \frac{-1-\sqrt{3}}{8}, \frac{3+\sqrt{3}}{8}\right),
\pm\left(\frac{-3-2 \sqrt{3}}{12}, \frac{4+\sqrt{3}}{12}\right)
\end{equation*}

\item Dodecagon: corona limit of $3.4.6.4;3^2.4.3.4$ (2-02)
\begin{equation*}
\pm\left(\frac{2+\sqrt{3}}{3}, 0  \right),
\pm\left(\frac{6+3\sqrt{3}}{10}, \frac{3+2 \sqrt{3}}{10}\right),
\pm\left(\frac{2+\sqrt{3}}{6}  , \frac{3+2 \sqrt{3}}{6} \right)
\end{equation*}
\begin{equation*}
\pm\left( 0, \frac{3+2 \sqrt{3}}{5}  \right),
\pm\left(\frac{-2-\sqrt{3}}{6}  , \frac{3+2 \sqrt{3}}{6}  \right),
\pm\left(\frac{-6-3\sqrt{3}}{10}, \frac{3+2 \sqrt{3}}{10} \right)
\end{equation*}

\item Hexadecagon: edge corona limit of $3^3.4^2;3^2.4.3.4$ (2-15)
\begin{equation*}
\pm\left( \frac{3+\sqrt{3}}{7} , 0 \right),
\pm\left( \frac{9+3\sqrt{3}}{22} , \frac{3+\sqrt{3}}{22} \right),
\pm\left( \frac{3+\sqrt{3}}{10} , \frac{3+\sqrt{3}}{10} \right),
\end{equation*}
\begin{equation*}
\pm\left( \frac{3+\sqrt{3}}{22} , \frac{9+3\sqrt{3}}{22} \right),
\pm\left( 0 , \frac{3+\sqrt{3}}{7}  \right),
\pm\left( \frac{-3-\sqrt{3}}{22}  , \frac{9+3\sqrt{3}}{22}  \right),
\end{equation*}
\begin{equation*}
\pm\left( \frac{-3-\sqrt{3}}{10} , \frac{3+\sqrt{3}}{10}  \right),
\pm\left( \frac{-9-3\sqrt{3}}{22}, \frac{3+\sqrt{3}}{22}\right).
\end{equation*}
\end{itemize}

\section{A repetitive tiling without a corona limit}
\label{NonUnif}
A translate of a patch $P=\{T_i\ |\ i\in F(P)\}$ by $x\in \R^l$ is
defined by $P+x=\{T_i+x\ |\ k\in F(P)\}$, where  
$P+x$ may not be a patch in $\T$. Two patches $P,Q$ are translationally
equivalent if there exists $x\in \R^l$ such that $Q=P+x$. 
A tiling $\T$ is {\it repetitive} if for any patch $P$, 
its translations appear infinitely often in $\T$. 
A tiling $\T$ has finite local complexity (FLC) if for any $R>0$, there are only finitely many 
patches in $B(x,R)$ with $x\in \R^l$ up to translation. 
Here is an example of a repetitive FLC Delone tiling which does not have a corona limit.

\begin{ex}
\label{RepNoLimit}
Consider a word monoid over a countable alphabet $\A=\{C_n\ |\ n=1,2,\dots\}$, whose binary operation is the concatenation of words. An empty word $\lambda$ is the identity. Define a monoid homomorphism $\sigma$ by:
$$
\sigma(C_n)=C_{n-1}C_nC_{n+1}C_nC_{n-1} \qquad n=1,2,\dots.
$$
Here we put $C_0=\lambda$. The action of $\sigma$ to a right infinite word $w_1w_2\dots$ is defined by $\sigma(w_1)\sigma(w_2)\dots$, and the same for the left infinite word. There is a fixed right infinite word 
$$
\sigma(w)=w=C_1C_2C_1C_1C_2C_3C_2C_1C_1C_2C_1C_1C_2C_1\dots
$$
and its mirror image is a left infinite word $w'$, which also satisfy $\sigma(w')=w'$.
We obtain a bi-infinite word $w'w\in \A^{\Z}$. 
Every finite subword in $w'w$ is a subword
of $\sigma^M(C_1)$ for some $M$, because the word $C_1C_1$ at the conjunction of $w'$ and $w$ is a subword of $\sigma^2(C_1)$. We can easily confirm by induction 
that $C_iC_j$ is a subword of $w'w$ then 
$i-j\in \{-1,0,1\}$, and $i-j=0$ happens only when $i=j=1$.
The word $C_1C_1C_1$ does not show. 
For each $C_i$ we associate an interval of length $2^{2^{i-1}}$. Prepare two 
intervals $A=[0,1]$ and $B=[0,2]$ and tile $C_i$ by $A$ if $i$ is odd, and by $B$
if $i$ is even. 
Tile the real line by intervals of length $2^{2^{n-1}}$ for $n=1,2,\dots$ 
according to the order of the word $w'w$, and then subdivide them by $A$, $B$ by this rule. Then we obtain a tiling of $\R$ by $A$ and $B$. If we see the word $BA^sB$ or $AB^tA$ in the final tiling, one can uniquely recover the word over $\A$ which produces $A^s$ or $B^t$. For e.g., $BB$ in $ABBA$ is produced by $C_2$, $AAAA$ in $BAAAAB$ is produced by $C_1C_1$.
This tiling is repetitive from
the above property of $w'w$. Since $2^{2^{n-1}}$ is rapid enough, 
it is easy to show that the directional speed does not exist. 
\end{ex}

\section{Problems and future perspectives}
\label{Final}
There are many intriguing problems in corona limits. Here we list down some of them: 

\begin{itemize} 
\item Is there a bound on the number of extreme points of corona limits for planar lattice periodic tilings?
Can we characterize the set of numbers of the extreme points? 
\item For which $n$ can we give a lattice periodic planar tiling whose corona limit is a regular $2n$-gon? {\it After all, we know that it happens only when $n=2,3$. See the paragraph after the proof of Theorem \ref{periodic-coronalimit} and \cite{ShuMal2014}.}
\item What can be said about corona limits of uniformly repetitive tilings? 
How about non-periodic self-similar tilings, or tilings generated by cut and projection? {\it As we discussed in the introduction, there are no universal method yet, but many partial results are known for concrete tilings.}
\item Is there a uniformly repetitive tiling whose corona limit is a ball?
\end{itemize}

Here a tiling is uniformly repetitive if
for any patch $P$, there exists a positive $r$ such that for any $x$, the ball $B(x,r)$ contains a translate of $P$. 

We finish this article by relating corona limits with the chemical and physical development of crystallization. 
Previous studies were devoted on the {\em structural forms} of crystals, where researchers look at the shape of crystals without regards to the external 
condition of the material.
Bravais \cite{Bravais1850} proposed a law which states that the faces most likely to be found on a crystal are those parallel to lattice planes of highest reticular density. The Bravais law relies only on the 32 point groups. Donnay and Harker \cite{DH1937} extended the rule to the 230 space groups. Hartman and Perdok~\cite{HP1955} focused on crystal zones and classified crystal faces according to the number of periodic bond chains (PBCs). 
It is remarkable that the structural form of a crystal is applicable also to the growth form, which we mention later. 
For instance, faces with two or more PBCs govern the growing shape of the crystal. However, since it is difficult to obtain bond chains in practice and to solve the problem, Bennema and Eerden \cite{BE1987} proposed the connected net model. 

Taking external conditions into account, there are two stages according to the macro-scale shapes of crystals. 
The first is about {\em equilibrium forms} and the second is about {\em growth forms} of crystals. 
Crystal growth needs to be regarded as the movement of a solid-liquid (gas) interface whose driving force is defined by the difference between solid and liquid (gas) chemical potentials at the interface per small distance. 

The equilibrium form of a crystal is defined as the final shape assumed by a growing crystal system as it arranges itself such that its surface Gibbs free energy is minimized. If its surface free energy density $\alpha(\theta)$ for each direction $\theta$ is known, the equilibrium form can be computed (it is known as the Wulff plot \cite{Wulff1901}). 
The growth form of a crystal is used to describe the growth process towards the equilibrium form. 
If the size of a crystal is small, its growth speed may depend on its surface tension. 
However, if its size is large enough, then the surface tension can be ignored. Moreover, if the transportation of atoms is fast enough, 
then the growth speed $dg/dt$ is proportional to its driving force, 
where the proportionality coefficient per atom depends on the direction $\theta$. We call this coefficient the
kinetic growth coefficient and denote this by $K(\theta)$. 
Under this assumption, Chernov \cite{Chernov1963} showed that the asymptotic form of a crystal is similarly enlarged and the form is independent of the shape of the initial nucleus. If the coefficient $K(\theta)$ is known, the asymptotic growth form can be computed in the same way as in the Wulff plot. The property of steady directional growth speed had been known as a good estimation since early times. For example, Kolmogorov \cite{Kolmogorov49} showed a  mathematical explanation of the mechanism of geometrical selections of crystals. But in the context of free energy minimization, the analysis of directional growth speed is difficult in general. 
For more details on the chemical aspects of crystals, see, for instance, \cite{Sunagawa}. 

Though shape study of aperiodic crystal seems not developed much, 
we find a few theoretical references. 
An analogy of Wulff-shapes, which optimizes the Gibbs-Curie surface energy, 
was defined and investigated in \cite{BorSch98,Sch2000}.
Shelling number and coordination number, meanwhile, are studied
in \cite{Baake-Grimm:97, Baake-Grimm:06}. These numbers give a certain average rate of growth. 

Corona limits in this paper are characterized using discrete patches (limit of the
$n$-th corona $P^{(n)}$). It is remarkable that as an analogy of the result by Chernov \cite{Chernov1963}, 
the shape of the corona limit is independent of the initial patch chosen (see Lemma \ref{Unique}). 
We introduced the notion of {\bf directional speed} $d(\v)$ 
in each direction $\v$, which is an analogy of the growth speed $dg/dt$. 
Therefore we expect that the corona limits $K$ can be regarded as models of crystal growth. 
We saw that the corona limit of a lattice periodic tiling is a convex polyhedron which is symmetric with respect to the origin (see Theorem \ref{periodic-coronalimit}). The shape of a corona limit depends on both the tiling and the adjacency relation defined on it. 
It is an interesting problem to reflect the growth condition of a crystal as its adjacency condition. 
The representation of directional speeds might help to bridge the theory of PBC and the traditional model based on growth coefficients.
We hope that this paper gives a new insight to the mathematical study of crystal shapes. 
\medskip

{\bf Acknowledgments.}
 We express our cordial gratitude to A.~V.~Shutov for informing us of the current status of the research on corona limits and for providing us some of the references, which are not easily accessible. We are also largely indebted to the anonymous referee and Fumihiko Nakano who gave us invaluable suggestions and related references.
 The first author is partially supported by JSPS grants (17K05159, 17H02849, BBD30028). The third author is partially supported by JSPS grants (26330016, 17K00015). The fourth author is supported by JSPS grant (15K17505).


\clearpage

\appendix
\section{1-uniform tilings and velocities}
\label{sec:1-regular-plots}

Computation of corona limits of $1$-uniform tilings using point adjacency and 
edge adjacency, see Remark \ref{adjacency} and Section \ref{Table}. Each row consists of 5 figures: tiling, (finite) coronas, their velocities, (finite) edge-coronas and their velocities. The convex hull of velocities is the corona limit, see Theorem \ref{periodic-coronalimit} and the description after it. 

\begin{figure}[h]
\includegraphics[scale=1]{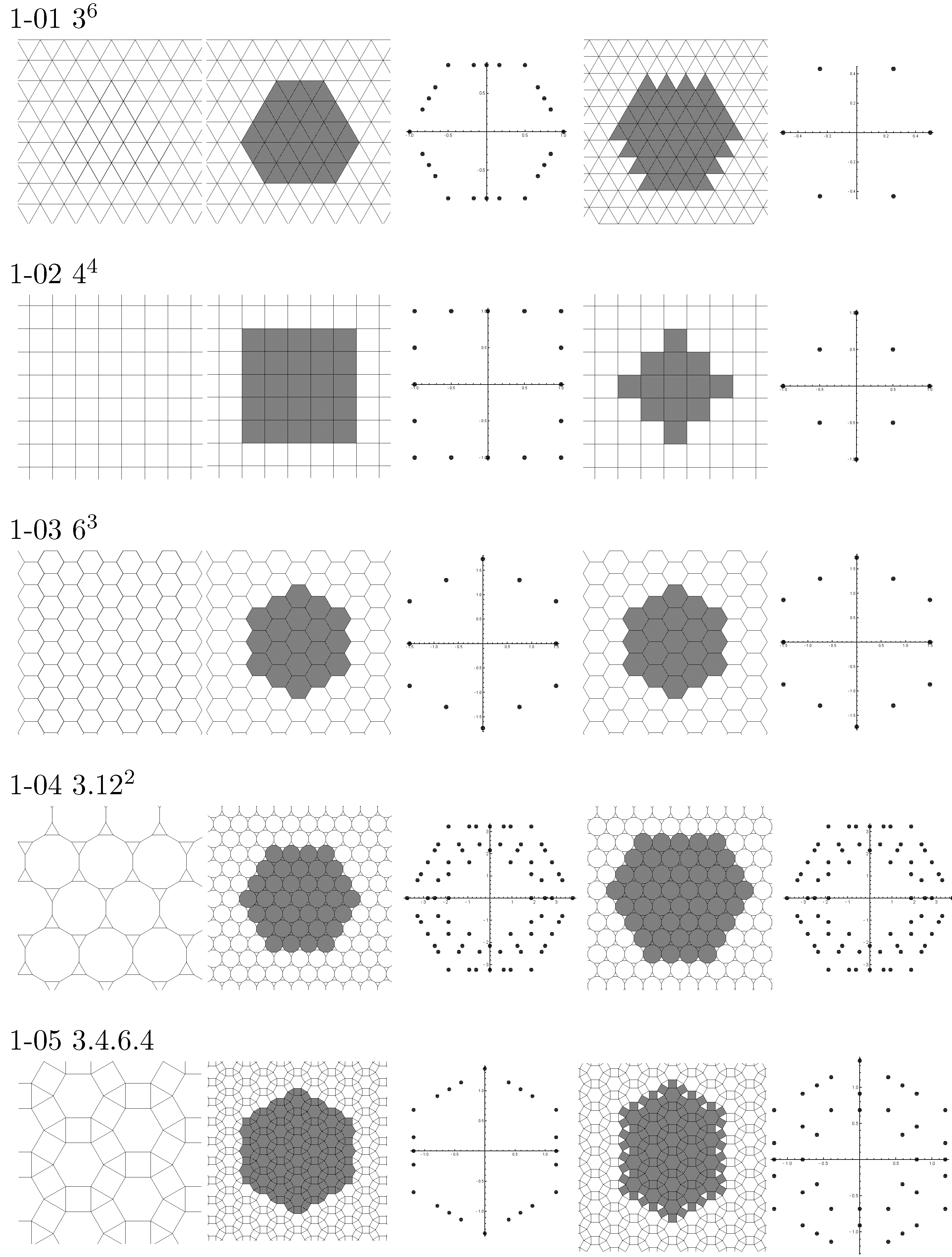}
\end{figure} 

\begin{figure}[h]
\includegraphics[scale=1]{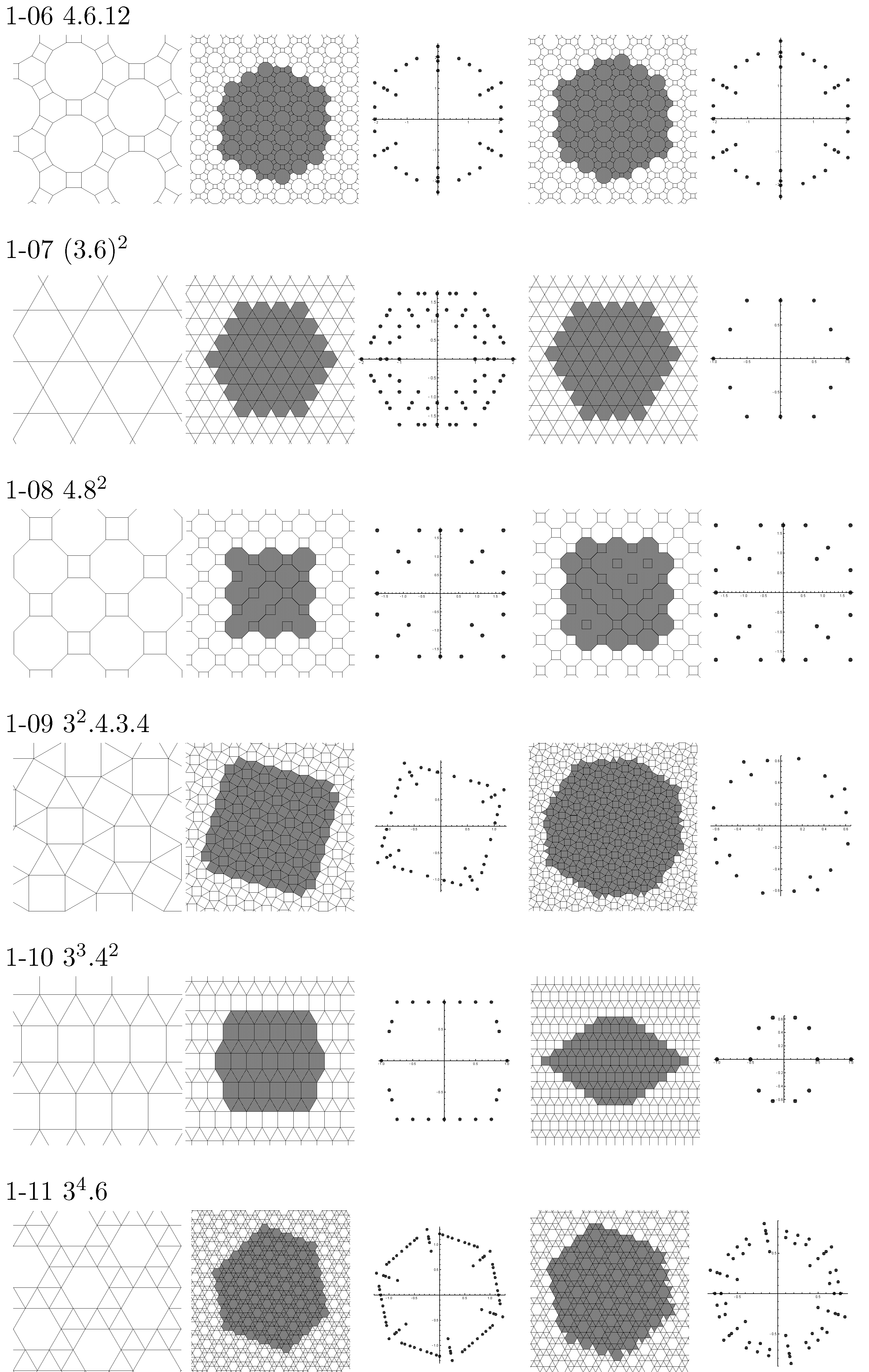}
\end{figure} 

\clearpage

\section{2-uniform tilings and velocities}
\label{sec:2-regular-plots}

The same computation of corona limits of $2$-uniform tilings. Each tiling is designated by two vertex configurations joined by semi-colon.

\begin{figure}[h]
\includegraphics[scale=1]{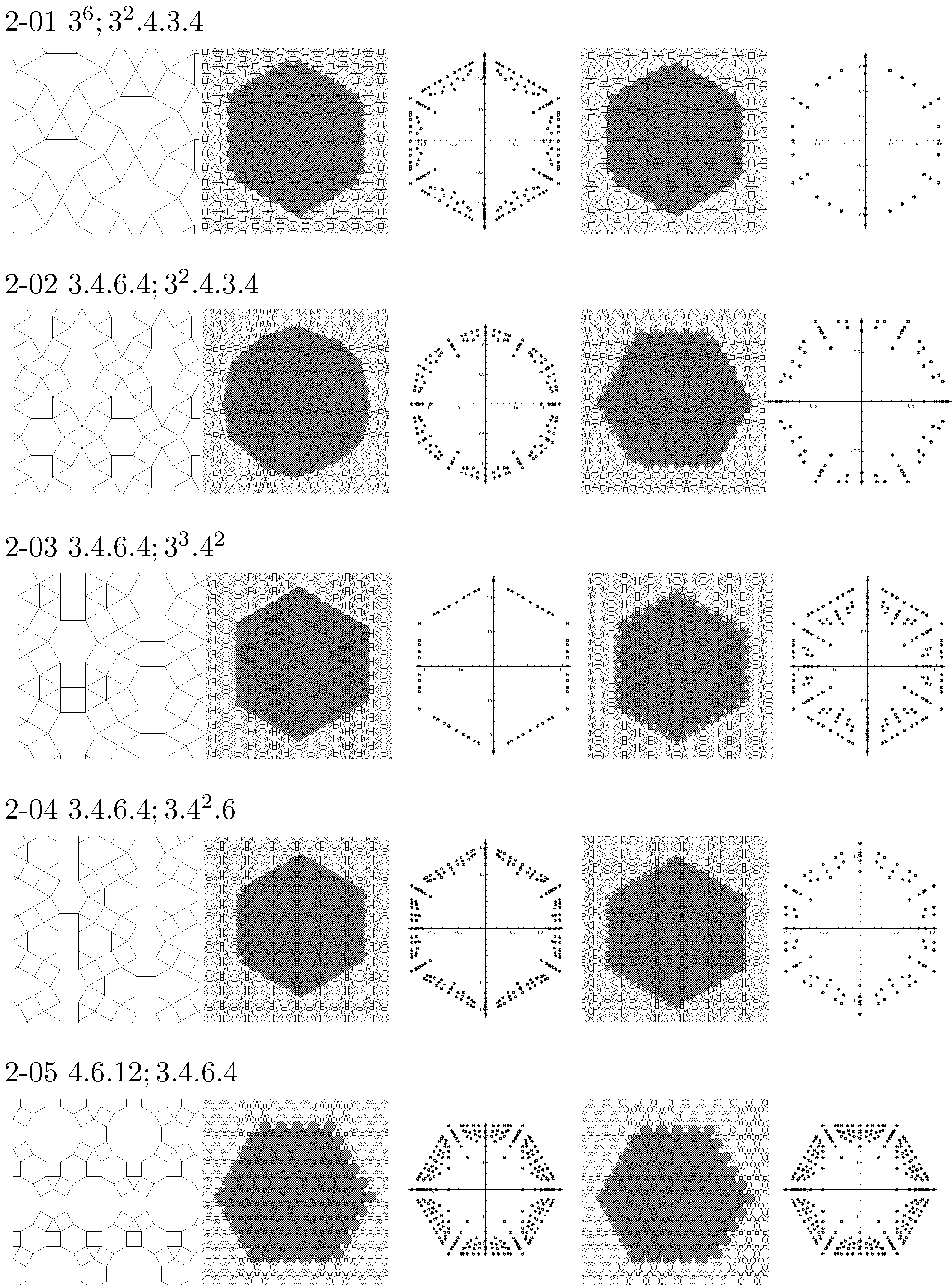}
\end{figure} 

\begin{figure}[h]
\includegraphics[scale=1]{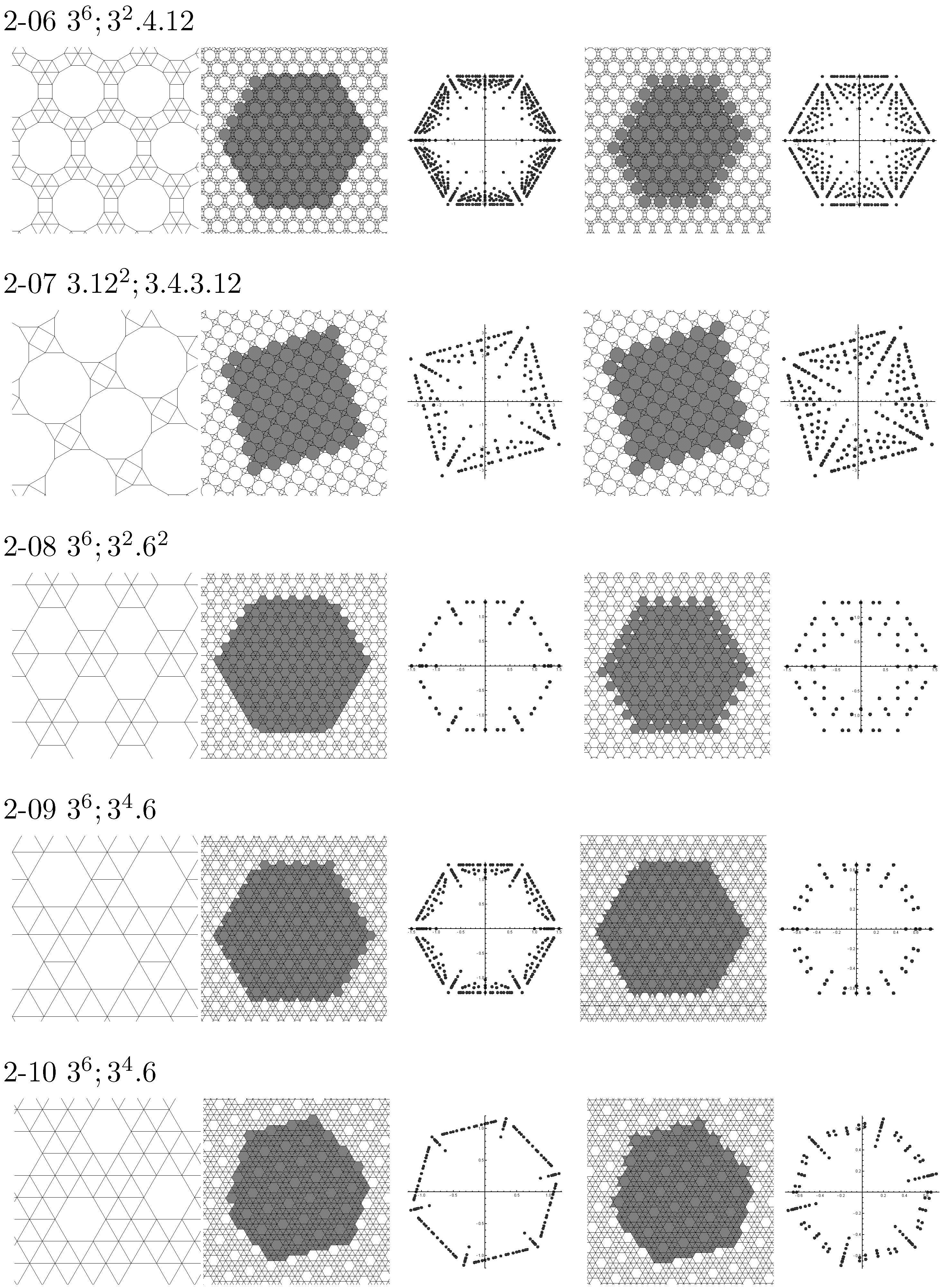}
\end{figure} 

\begin{figure}[h]
\includegraphics[scale=1]{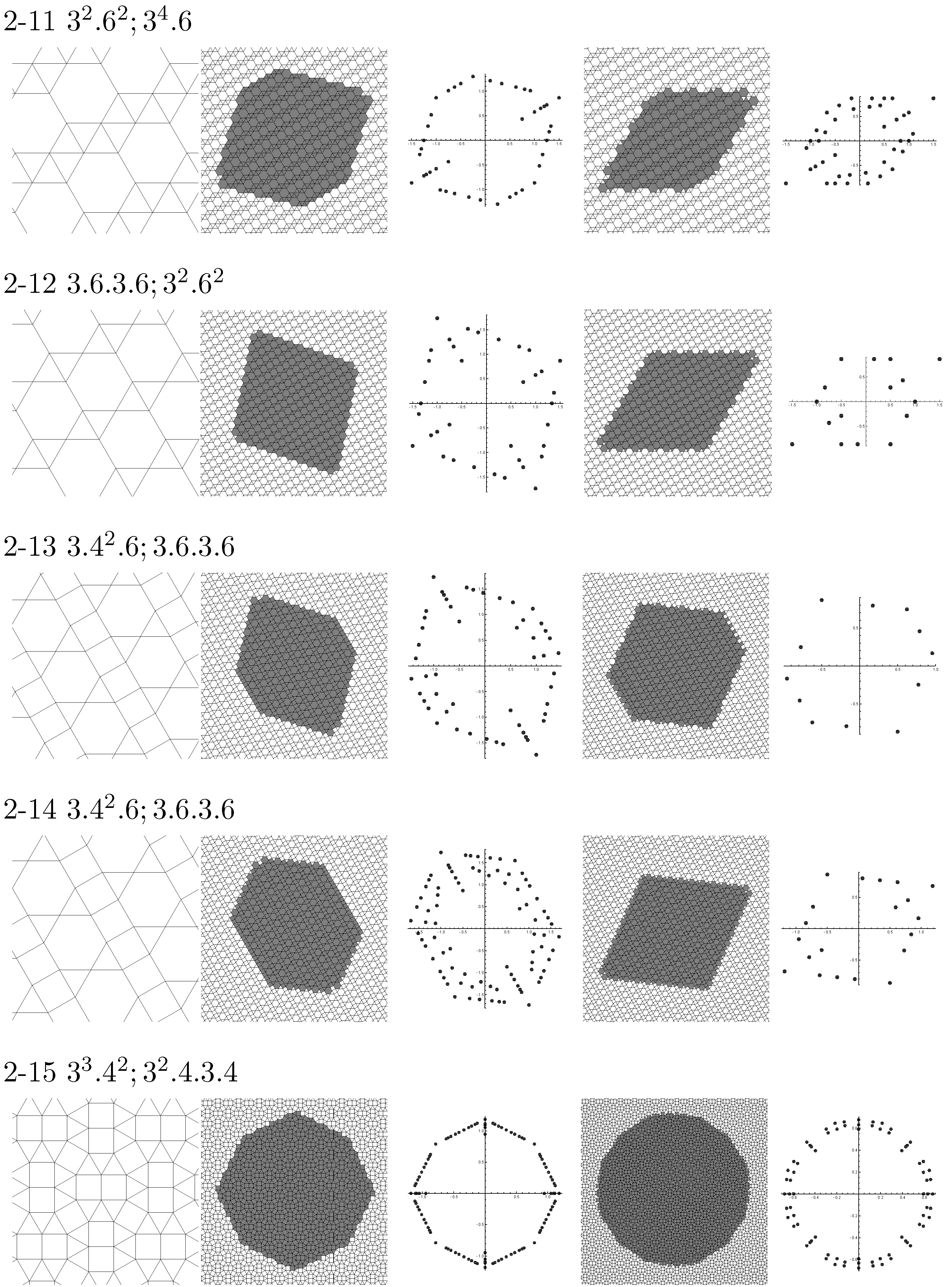}
\end{figure} 

\begin{figure}[h]
\includegraphics[scale=1]{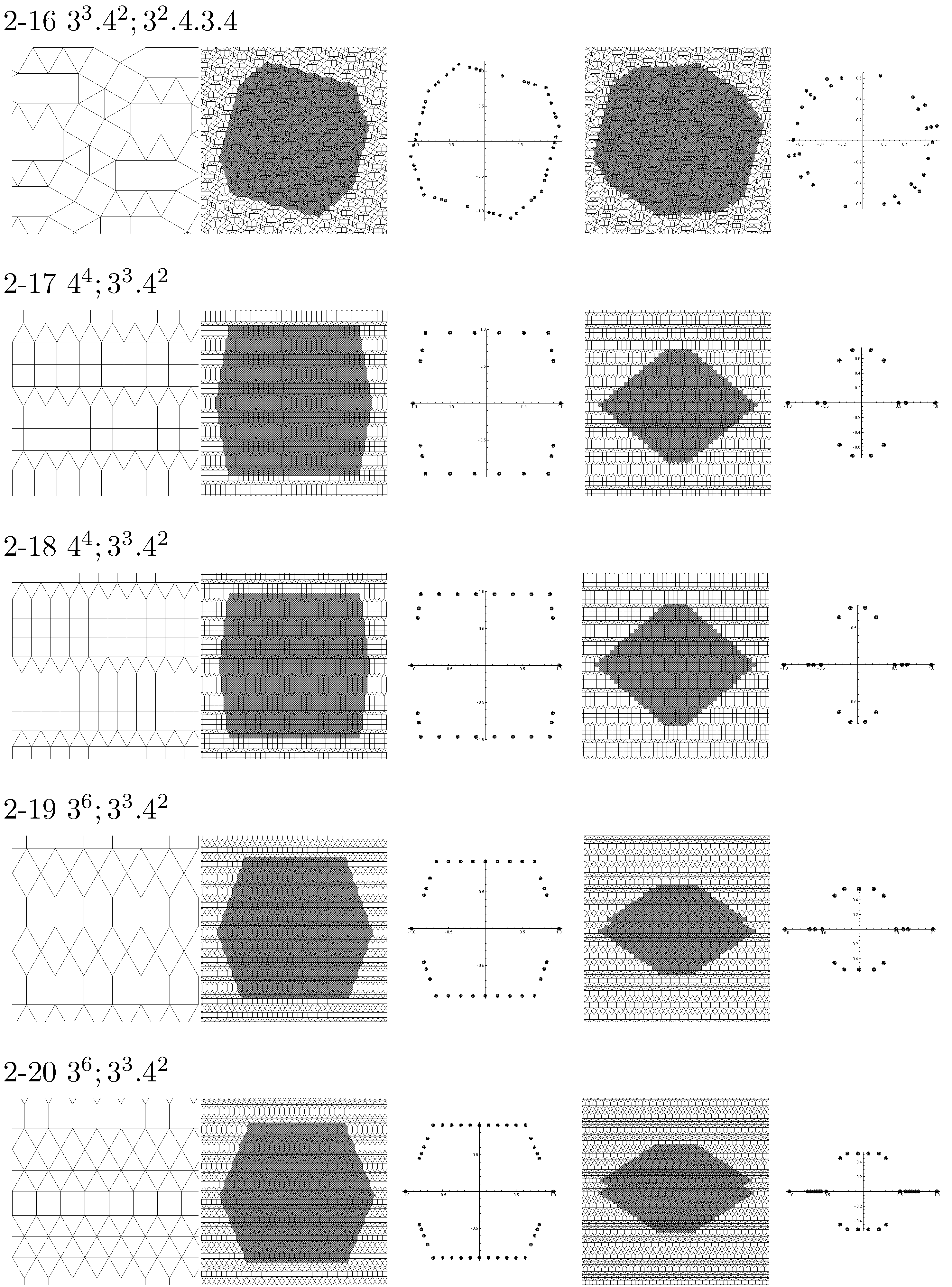}
\end{figure}

\end{document}